\documentclass[11pt,a4paper]{article}

\usepackage{amsmath}
\usepackage{amssymb}
\usepackage{graphicx}
\newtheorem{thm}{Theorem.}[section]

\usepackage{subfigure}

\author{Atsushi Iwasaki\thanks{Kyoto University\ \  e-mail: iwasaki.atsushi.4x@kyoto-u.ac.jp}}
\title{Independent Randomness Tests based on the Orthogonalized Non-overlapping Template Matching Test}
\begin{document}

\maketitle
\begin{abstract}
In general, randomness tests included in a test suite are not independent of each other.
This renders it difficult to fix a rational criterion through the whole test suite with an explicit significance level.
In this paper, we focus on the Non-overlapping Template Matching Test, which is a randomness test included in the NIST statistical test suite.
The test uses a parameter called ``template'' and we can consider a test item for each template.
We investigate dependency between two test items by deriving the joint probability density function of the two p-values and propose a transformation to make multi test items independent of each other.
\end{abstract}

\section{Introduction}
Random number sequences are widely used in many fields such as information security and various stochastic algorithms.
In general, the properties of the random number sequences used in these applications affects their efficiency and effectiveness.
In information security contexts including cryptography, the quality of randomness is particularly important because security strength depends on this quality.
Thus, it is necessary and important to evaluate the randomness of sequences or their generators from different perspectives.

Randomness tests are one such means of evaluation.
They are hypothesis tests and the null is that the given sequence is truly random.
Randomness tests do not require information about the generator of the given sequence.
Thus, they can be widely used without regard to the generator.
Recently, Tamura and Shikano used randomness tests to inspect the properties of a quantum computer developed by IBM and showed that this computer does not work as expected \cite{Tamura-Shikano}.

There are many types of randomness tests and particular test suites have been proposed \cite{Diehard, TestU01, NIST}.
SP800-22 \cite{NIST} published by the US National Institute of Standard and Technology (NIST) is one of the most well-known test suites.
The first version was published in 2001 and revision 1a published in 2010 is currently valid.
Revision 1a consists of 188 test items in the default setup that can be categorized into 15 test types.
Some of the tests included therein are parametric, and so multi items can be implemented for each test.

For each test item included in SP800-22, a criterion is specified concerning whether the given sequences pass or not.
However, a criterion that simultaneously covers all test items is not specified.
This is a critical problem common to many test suites that needs to be addressed for the effective use of randomness tests. 
In general, test items included in one test suite are not independent of each other, i.e., p-values computed by multi test items do not distribute independently and the joint distribution is not known.
This renders it difficult to specify a rational criterion through all test items.
As an exception, Sugita proposed a test suite consisting of items that are independent of each other \cite{Sugita}, where p-values are computed based on ``parity''.
Previous studies have reported empirical results concerning the dependency among particular  randomness tests \cite{Hellekalek,Turan,DO,Fan,Yamaguchi2016,Sulak,DO2017,Iwasaki2018}. These studies were conducted under the null or under a certain alternative hypothesis. Considering  dependency under an alternative hypothesis can be important from the perspective of computational burden.
Notwithstanding the importance of extant work in this domain, the dependency among test items included in a given test suite has not received sufficient attention in the literature.
In this paper, we aim to address this by focusing on test suite dependency under the null, culminating in the specification of a criterion through all test items.

To specify a criterion through all test items in SP800-22, it is critically important to study the Non-overlapping Template Matching Test, which is one type of randomness test included in SP800-22.
This test counts the number of occurrences of a short string called ``template'' on a given sequence and computes a p-value based on the number.
Various short strings can be used as templates, and so we can consider multi test items in the context of the Non-overlapping Template Matching Test.
Such test items account for 148 of the 188 items in SP800-22. Thus, understanding the dependency among these 148 test items is fundamental for understanding the dependency among test items in SP800-22.
We investigate the dependency between two Non-overlapping Template Matching Test items by deriving the joint probability density function of the two p-values.
Next, we propose a transformation to render multi test items independent of each other.

The remainder of the paper is organized as follows.
In section 2, we introduce the Non-overlapping Template Matching Test and discuss its soundness.
In section 3, we derive the joint distribution of two Non-overlapping Template Matching Test p-values.
In section 4, we propose a method to render the test items independent.
Finally, in section 5 we offer conclusions.

\section{Non-overlapping Template Matching Test}

In this section, we introduce the Non-overlapping Template Matching Test and confirm its soundness.
Here, ``soundness'' means that a p-value follows the uniform distribution on $[0,1]$ under the null hypothesis, at least approximately.

\subsection{Algorithm}
The Non-overlapping Template Matching Test algorithm is as follows:
\begin{enumerate}
\item Divide the given $n$-bit sequence into $N$ blocks $B^{(1)}, B^{(2)},\cdots,B^{(N)}$.
Here, each block is $M$-bit and $M=\left\lfloor\frac{n}{N}\right\rfloor$.
\item For $j=1,2,\cdots,N$, compute $c_j$ as the occurrences of a given short string $T$ called ``template" on $B^{(j)}$, i.e.,
\begin{align}
\label{siki2-0}
c_j=\#\left\{k\in\{1,2,\cdots,M-m+1\}\mid  {B^{(j)}}_{[k,k+m-1]}=T\right\}.
\end{align}
Here, $X_{[s,t]}$ means the $(t-s+1)$-bit substring from $s$-th bit to $t$-th bit of a string $X$. 
\item Compute $\chi_{obs}$ as follows:
\begin{align}
\chi_{obs}=\sum_{j=1}^N\left(\frac{c_j-\mu}{\sigma}\right)^2,
\end{align}
where
\begin{align}
\mu=&\frac{M-m+1}{2^m},\\
\sigma=&M\left(\frac{1}{2^m}-\frac{2m-1}{2^{2m}}\right).
\end{align}
\item Perform a $\chi^2$-test with $N$ degrees of freedom for $\chi_{obs}$ and compute p-value $p$.
\end{enumerate}
We can use any short string $T$ as a template provided $T$ satisfies 
\begin{align}
\label{siki2-1}
T_{[1+k,m]}\ne T_{[1,m-k]}
\end{align}
for $k=1,2,3,\cdots,m-1$.
Where $m=9$, there are 148 templates satisfying the condition.
 The sample program provided by NIST uses the 148 9-bit templates in the default set up.
The program performs these 148 test items as Non-overlapping Template Matching Tests, i.e., each test item is operationalized via the above algorithm using one of the 148 templates.

\subsection{Soundness}
To satisfy the property of soundness, $\chi_{obs}$ is required to at least approximately follow a $\chi^2$-distribution with $N$ degrees of freedom, under the null hypothesis.
The requirement is satisfied if each $c_j$ independently follows a normal distribution $\mathcal{N}(\mu,\sigma^2)$, at least approximately.
Since the blocks do not overlap with each other, independence is assured.
Then, we show that $c_j$ follows $\mathcal{N}(\mu,\sigma^2)$.
 This also serves as preparation for the next section. 

First, we confirm that the average and the variance of $c_j$ are $\mu$ and $\sigma^2$, respectively.
We regard $B^{(j)}$ as a $M$-bit random variable uniformly distributed on $\{0,1\}^M$ and let $I_k$ be a random variable described as
\begin{align}
\label{siki2-2}
I_k=\begin{cases}
1\ \ &({B^{(j)}}_{[k,k+m-1]}=T)\\
0\ \ &({B^{(j)}}_{[k,k+m-1]}\ne T)
\end{cases}
.
\end{align}
Since $B^{(j)}$ is uniformly distributed, it is straightforward that
\begin{align}
\text{Prob}\{I_k=1\}=&\frac{1}{2^m},\label{siki2-2-1}\\
\text{Prob}\{I_k=0\}=&1-\frac{1}{2^m}\label{siki2-2-2}.
\end{align}
Using $\{I_k\}_{k=1}^{M-m+1}$, we can represent $c_j$ as
\begin{align}
\label{siki2-0}
c_j=\sum_{k=1}^{M-m+1}I_k.
\end{align}
Then, by (\ref{siki2-2-1}) and (\ref{siki2-2-2}),
\begin{align}
\mathbb{E}\left[c_j\right]=&\mathbb{E}\left[\sum_{k=1}^{M-m+1}I_k\right]\\
=&\sum_{k=1}^{M-m+1}\mathbb{E}\left[I_k\right]\\
=&\sum_{k=1}^{M-m+1}\frac{1}{2^m}\\
=&\frac{M-m+1}{2^m}\\
=&\mu.
\end{align}
Similarly, the standard deviation can be derived as follows:
\begin{align}
\mathbb{E}\left[c_j^2\right]-{\mathbb{E}\left[c_j\right]}^2
=&\mathbb{E}\left[\sum_{k=1}^{M-m+1}\sum_{l=1}^{M-m+1}I_kI_l\right]-\left(\frac{M-m+1}{2^m}\right)^2\\
=&\sum_{k=1}^{M-m+1}\sum_{l=1}^{M-m+1}\mathbb{E}\left[I_kI_l\right]-\left(\frac{M-m+1}{2^m}\right)^2\label{siki2-3}.
\end{align}
By (\ref{siki2-1}) and (\ref{siki2-2}),
\begin{align}
\label{siki2-4}
\mathbb{E}\left[I_kI_l\right]=
\begin{cases}
\frac{1}{2^m}\ \ &(k=l)\\
0\ \ &(1\leq |k-l|\leq m-1)\\
\frac{1}{2^{2m}}\ \ &(|k-l|\geq m)
\end{cases}
.
\end{align}
Substituting (\ref{siki2-4}) into (\ref{siki2-3}), we obtain
\begin{align}
&\mathbb{E}\left[c_j^2\right]-{\mathbb{E}\left[c_j\right]}^2\\
=&\left\{\frac{M}{2^m}-\frac{(M-m+1)(M-3m+2)}{2^{2m}}+O\left(\frac{m^2}{2^m}\right)\right\}-\left(\frac{M-m+1}{2^m}\right)^2\\
=&M\left(\frac{1}{2^m}-\frac{2m-1}{2^{2m}}\right)+O\left(\frac{m^2}{2^m}\right).\label{sigma}
\end{align}
Then, $\sigma^2$ approximates the variance of $c_j$ for sufficiently large $M$.

Next, we confirm that $c_j$ approximately follows a normal distribution.
We use the following theorem \cite{Hoeffding}:
\begin{thm}
\label{thm1}
Consider a random variable sequence $\{d_k\}_{k=1}^\infty$ satisfying the following conditions:
\begin{itemize}
\item There exists a positive integer $s$ such that $\{d_k\}_{k=1}^t$ and $\{d_k\}_{k=s+t}^\infty$ are independent of each other for all $t$.
\item For an arbitrary positive integer $u$, the joint distribution of $\{d_k\}_{k=v}^{v+u}$ does not depend on $v$.
\item  $\mathbb{E}\left[|d_1|^3\right]<\infty$.
\end{itemize}
Then, the distribution of $D_a:=\sum_{k=1}^a\left(d_k-\mathbb{E}\left[d_k\right]\right)/\sqrt{a}$ converges to a normal distribution as $a\to\infty$.
\end{thm}
As $M\to\infty$, the random variable sequence $\{I_k\}_{k=1}^{M-m+1}$ satisfies the three conditions in theorem \ref{thm1}.
Then, by (\ref{siki2-0}) it is shown that $c_j$ approximately follows a normal distribution.

From the above, it is proven that $c_j$ follows $\mathcal{N}(\mu,\sigma^2)$ for sufficiently large $M$.
Consequently, the soundness of the Non-overlapping Template Matching Test is established.

\section{Joint distribution of p-values}
Consider two $m$-bit templates $T^{(1)}$ and $T^{(2)}$. 
In the following, we use the notation $X^{(s)}$ to denote variable $X$ corresponding to $T^{(s)}$.
The purpose of this section is to obtain the joint distribution of $(p^{(1)},p^{(2)})$ under the null hypothesis.

\subsection{Derivation of the joint distribution}
First, we derive the joint distribution of $\left(\frac{c_j^{(1)}-\mu}{\sigma},\frac{c_j^{(2)}-\mu}{\sigma}\right)$.
Theorem \ref{thm1} can be extended to multi-dimensional cases and thus it is shown that$\left(\frac{c_j^{(1)}-\mu}{\sigma},\frac{c_j^{(2)}-\mu}{\sigma}\right)$ follows a two-dimensional normal distribution as $M\to\infty$.
Since we already know that both marginal distributions follow a standard normal distribution, the joint distribution is entirely specified if we derive the correlation coefficient.
Using $I_k^{(1)}$ and $I_k^{(2)}$, we obtain
\begin{align}
&\mathbb{E}\left[c_j^{(1)}c_j^{(2)}\right]\\
=&\mathbb{E}\left[\sum_{k=1}^{M-m+1}\sum_{l=1}^{M-m+1}I_k^{(1)}I_l^{(2)}\right]\\
=&\sum_{k=1}^{M-m+1}\sum_{l=1}^{M-m+1}\mathbb{E}\left[I_k^{(1)}I_l^{(2)}\right].\label{3-1}
\end{align}
We define the following notation:
\begin{align}
e_k:=&\begin{cases}
1\ \ &(1\leq k\leq m-1\text{ and }{T^{(1)}}_{[1,m-k]}={T^{(2)}}_{[1+k,m]})\\
1\ \ &(-m+1\leq k\leq -1\text{ and }{T^{(1)}}_{[1-k,m]}={T^{(2)}}_{[1,m+k]})\\
0\ \ &(\text{ otherwise })
\end{cases}
.
\end{align}
Then, we obtain
\begin{align}
\label{3-2}
\mathbb{E}\left[I_k^{(1)}I_l^{(2)}\right]=\begin{cases}
0\ \ &(k=l)\\
\frac{e_{k-l}}{2^{m+|k-l|}}\ \ &(1\leq |k-l|\leq m-1)\\
\frac{1}{2^{2m}}\ \ &(|k-l|\geq m)
\end{cases}
.
\end{align}
Substituting (\ref{3-2}) into (\ref{3-1}), we have
\begin{align}
&\mathbb{E}\left[c_j^{(1)}c_j^{(2)}\right]\\
=&\sum_{k=1}^{M-m+1}\mathbb{E}\left[I_k^{(1)}I_k^{(2)}\right]+\sum_{k=1}^{M-m+1}\sum_{l:\ 1\leq|k-l|\leq m-1}\mathbb{E}\left[I_k^{(1)}I_l^{(2)}\right]+\sum_{k=1}^{M-m+1}\sum_{l:\ |k-l|\geq m}\mathbb{E}\left[I_k^{(1)}I_l^{(2)}\right]\\
=&M\sum_{k=1}^{m-1}\frac{e_k+e_{-k}}{2^{m+k}}+\frac{(M-m+1)(M-3m+2)}{2^{2m}}+O\left(\frac{m^2}{2^m}\right).
\end{align}
From the above, the correlation coefficient of the joint distribution $\rho\left(T^{(1)},T^{(2)}\right)$ is computed as follows:
\begin{align}
\rho\left(T^{(1)},T^{(2)}\right)
=&\mathbb{E}\left[\frac{c_j^{(1)}-\mu}{\sigma}\cdot\frac{c_j^{(2)}-\mu}{\sigma}\right]\\
=&\frac{\mathbb{E}\left[c_j^{(1)}c_j^{(2)}\right]-\mu^2}{\sigma^2}\\
=&\frac{-2m+1+\sum_{k=1}^{m-1}2^{m-k}(e_k+e_{-k})}{2^m-2m+1}+O\left(\frac{m^2}{2^mM}\right).\label{siki3-3}
\end{align}
Hereafter, we ignore the remainder term of (\ref{siki3-3}).

Next, we consider the joint cumulative distribution of $\left(\chi_{obs}^{(1)},\chi_{obs}^{(2)}\right)$.
Assume that random variable pairs $(x_1,y_1)$, $(x_2,y_2)$, $\cdots$, $(x_N,y_N)$ independently follow a two-dimensional normal distribution and that each marginal distribution is a standard normal distribution.
We denote the correlation coefficient between $x_j$ and $y_j$ as $\rho$, which does not depend on $j$.
We define $X_N$ and $Y_N$ as follows:
\begin{align}
X_N=\sum_{j=1}^Nx_j^2,\\
Y_N=\sum_{j=1}^Ny_j^2.
\end{align}
It is straightforward that
\begin{align}
\text{Pr}\left\{X_N\leq0\right\}=\text{Pr}\left\{Y_N\leq0\right\}=0.
\end{align}
Thus, we consider only the case that $X_N>0$ and $Y_N>0$.
In the sample program provided by NIST, the number of blocks is fixed at 8, and so we adopt the assumption that $N$ is even.
When $-1<\rho<1$, the joint cumulative distribution function $F_{N,\rho}$ is as follows:
\begin{align}
F_{N,\rho}(X_N,Y_N)=\sum_{r=0}^\infty\frac{\rho^{2r}(1-\rho^2)^\frac{N}{2}}{r!\Gamma\left(r+\frac{N}{2}\right)\Gamma\left(\frac{N}{2}\right)}\gamma\left(r+\frac{N}{2},\frac{X_N}{2(1-\rho^2)}\right)\gamma\left(r+\frac{N}{2},\frac{Y_N}{2(1-\rho^2)}\right),
\label{siki3-4}
\end{align}
where $\gamma$ is the incomplete gamma function. 
Derivation of (\ref{siki3-4}) is elaborated in the Appendix.
Using  $F_{N,\rho}$, the joint cumulative distribution of  $\left(\chi_{obs}^{(1)},\chi_{obs}^{(2)}\right)$ is given as $F_{N,\rho\left(T^{(1)},T^{(2)}\right)}\left(\chi_{obs}^{(1)},\chi_{obs}^{(2)}\right)$.

Finally, we obtain the joint distribution of $\left(p^{(1)},p^{(2)}\right)$ using the joint cumulative distribution of  $\left(\chi_{obs}^{(1)},\chi_{obs}^{(2)}\right)$.
Assume that $\chi^2(P)$ is the value of $\chi_{obs}$ when the value of the corresponding p-value is $P$.
By the definition of the p-value, we have
\begin{align}
\text{Prob}\left\{p^{(1)}>P\ \text{and}\ p^{(2)}>P^\prime\right\}=F_{N,\rho\left(T^{(1)},T^{(2)}\right)}\left(\chi^2(P),\chi^2(P^\prime)\right).
\label{siki3-5}
\end{align}

\subsection{Experiment}
We generated $10^6$ sequences using the Mersenne twister \cite{MT}.
The length of each sequence is $10^6$-bit.
Using the sequences we computed two-dimensional joint distributions of p-values and compared the experimental and theoretical distributions based on (\ref{siki3-5}).
We used templates pairs $(001010101,010101011)$ and $(001010101,101010100)$ for expository purposes.
The correlation coefficients defined in the previous subsection can be computed as $\rho(001010101,010101011)=0.652525
$ and $\rho(001010101,101010100)=0.321212$.

\begin{figure}[thbp]
\begin{center}
\subfigure[Experiment]{
\includegraphics[width=75mm]
{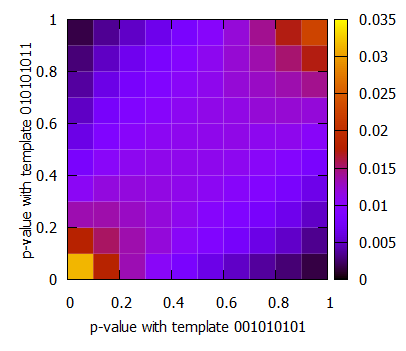}
\label{left}}
\subfigure[Theory]{
\includegraphics[width=75mm]
{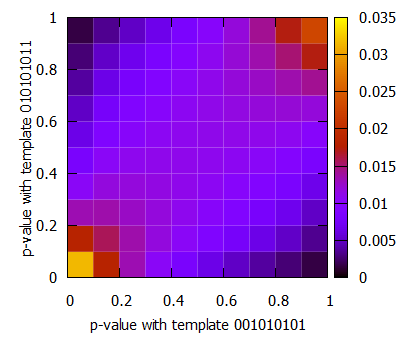}
\label{righy}}
\caption{Joint distribution of p-values with templates 001010101 and 010101011}
\label{fig1}
\end{center}
\end{figure}
\begin{figure}[thbp]
\begin{center}
\subfigure[Experiment]{
\includegraphics[width=75mm]
{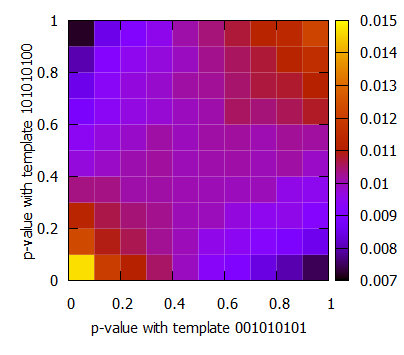}
\label{left}}
\subfigure[Theory]{
\includegraphics[width=75mm]
{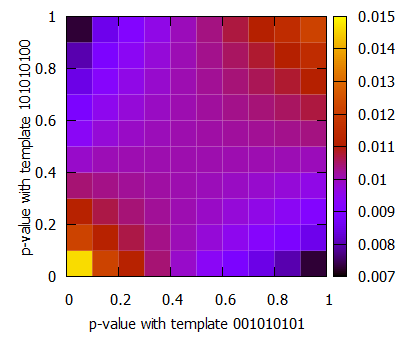}
\label{righy}}
\caption{Joint distribution of p-values with templates 001010101 and 101010100}
\label{fig2}
\end{center}
\end{figure}

The results are illustrated in Figure \ref{fig1} and \ref{fig2} and testify to the robustness of the theory discussed in this section. 
It is also shown that we cannot regard the test items as independent of each other in practical contexts.

\section{Orthogonalization}

For each sequence, the sample program provided by NIST uses 148 templates and computes 148 Non-overlapping Template Matching Test p-values in the default setup.
Thus, we need to extend the discussion in the former section and retrieve the 148 dimensional joint distribution of the 148 p-values.
However, 148 dimensions is computationally onerous, so this needs to be circumvented.
To transform the 148 test items and generate new test items that are independent of each other will be one of the solutions.

\subsection{Proposed method}
We propose a method to transform the test items by orthogonalization of a multi-dimensional normal distribution.
We consider $m$-bit templates $T^{(1)}$, $T^{(2)}$, $\cdots$, $T^{(R)}$ and introduce notation ${\mathbf C}_j$ as
\begin{align*}
{\mathbf C}_j:=\left(\frac{c_j^{(1)}-\mu}{\sigma},\frac{c_j^{(2)}-\mu}{\sigma},\cdots,\frac{c_j^{(R)}-\mu}{\sigma}\right)^\top.
\end{align*}
Assume that $M$ is sufficiently large.
The covariance matrix of ${\mathbf C}_j$ denoted by $\Sigma$ is described as follows:
\begin{align*}
(\Sigma)_{k,l}=
\begin{cases}
1\ \ (k=l)\\
\rho(T^{(k)},T^{(l)})\ \ (k\ne l)
\end{cases}
.
\end{align*}
If $\text{det}(\Sigma)\ne 0$, then the probability density function of ${\mathbf C}_j$ is described as
\begin{align*}
p({\mathbf C}_j)=\frac{1}{(2\pi)^\frac{R}{2}\sqrt{|\text{det}(\Sigma)|}}\exp\left(-\frac{{{\mathbf C}_j}^\top\Sigma^{-1}{{\mathbf C}_j}}{2}\right)
\end{align*}
when $M\to\infty$.
Since $\Sigma$ is a real symmetric matrix, then there exists an orthogonal matrix $L$ such that
\begin{align}
L^\top\Sigma^{-1} L=\Lambda
\end{align}
where $\Lambda$ is a diagonal matrix.
If $\text{det}(\Sigma)\ne 0$, then the diagonal elements of $\Lambda$ are all positive because $\Sigma$ is the covariance matrix of ${\mathbf C}_j$.
Then, using a diagonal matrix $Q$ described as
\begin{align}
\forall k,\ \ \ \ (Q)_{k.k}:=\frac{1}{\sqrt{(\Lambda)_{k,k}}},
\end{align}
we get 
\begin{align}
Q^\top L^\top\Sigma^{-1} LQ=I
\end{align}
where $I$ is a unit matrix.
Using $L$ and $Q$, we transform ${\mathbf C}_j$ to ${\mathbf C}^\prime_j$ as 
\begin{align}
{\mathbf C}^\prime_j:=(LQ)^{-1}{\mathbf C}_j.
\end{align}
Then, we obtain the probability density function of ${\mathbf C}^\prime_j$ as follows:
\begin{align}
p^\prime({\mathbf C}^\prime_j)=&p(LQ{\mathbf C}^\prime_j)\times |\text{det}(LQ)|\\
=&p(LQ{\mathbf C}^\prime_j)\times \sqrt{|\text{det}(\Sigma)|}\\
=&\frac{1}{(2\pi)^\frac{R}{2}}\exp\left(-\frac{{{\mathbf C}^\prime_j}^\top{{\mathbf C}^\prime_j}}{2}\right).
\end{align}
Then, components of ${\mathbf C}^\prime_j$ follow the standard normal distribution independent of each other.
For each $j$, if we replace ${\mathbf C}_j$ with ${\mathbf C}^\prime_j$, then the test items become independent of each other. 

However, it is not ensured that $\text{det}(\Sigma)\ne 0$.
Indeed, $\Sigma$ defined on 148 9-bit templates has zero eigenvalues.
Consider templates $T^{(1)}=100000000$ and $T^{(2)} =000000001$.
If ${B^{(j)}}_{[k,k+8]}=T^{(1)}$, then it is ensured that an integer $l>k+8$ exists such that
\begin{align}
{B^{(j)}}_{[k,l]}=10000\cdots001\ \ \text{, i.e.,}\ \ {B^{(j)}}_{[l-8,l]}=T^{(2)}
\end{align}
as $M\to\infty$.
In other words, where $T^{(1)}$ and $T^{(2)}$ occur they form a pair.
Then, we have
\begin{align}
\label{kousoku}
\lim_{M\to\infty}\frac{c_j^{(1)}-\mu}{\sigma}=\lim_{M\to\infty}\frac{c_j^{(2)}-\mu}{\sigma}.
\end{align}
Equation (\ref{kousoku}) is one reason why $\Sigma$ has zero eigenvalues.
Thus, we need to remove at least $T^{(1)}$ or $T^{(2)}$ from our discussion. 
Technically, we can identify templates to remove by checking eigenvectors corresponding to zero eigenvalues.
When we consider the 148 9-bit templates, we need to remove either $100000000$ or $000000001$, and $011111111$ or $111111110$. Finally, we need to remove $001010101$ or $010101011$ or $101010100$ or $110101010$.

\subsection{Experiment}
We used the Mersenne twister \cite{MT} and AES-128 \cite{AES} and generated $10^6$ sequences for each generator.
The length of the sequences is $10^6$-bit. 
AES-128 was used with counter mode.
We executed 145 test items for these sequences using 145 9-bit templates (except 100000000, 111111110 and 001010101) and compared the pre- and post-transformation results.

\begin{figure}
\begin{center}
\includegraphics{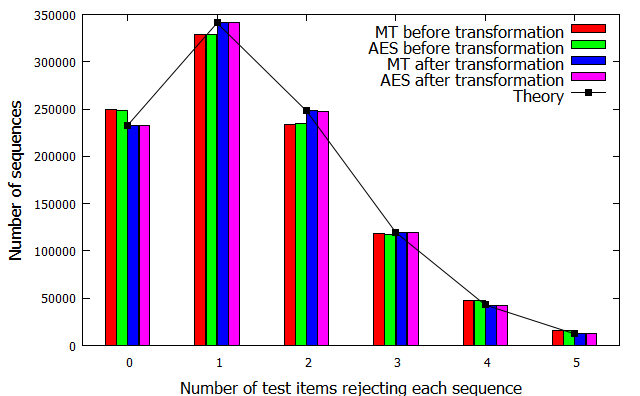}
\caption{Number of rejected sequences before and after the transformation}
\label{fig3}
\end{center}
\end{figure}

Figure \ref{fig3} shows the frequency of the sequences for the number of test items rejecting each sequence.
Here, ``reject'' means that a p-value is less than 0.01.
The black line denotes expectation values in the case that all test items are independent of each other.
This result implies that we can retrieve test item series that are independent of each other, and thus the proposed transformation is effective.

\section{Conclusions}
We derived the joint distribution of two Non-overlapping Template Matching Test p-values under the assumption that the block size is infinity.
The results suggest that we cannot regard test items as independent of each other.

We also proposed a method to remove dependency using orthogonalization of a multi-dimensional normal distribution.
Experimental results testify to the efficacy of the method.
Thus, it is expected that the method will contribute to fixing a rational criterion through all test items in SP800-22 and thus to the appropriate use of randomness tests.

\setcounter{section}{0}
\renewcommand{\appendix}{}
\appendix
\renewcommand{\thesection}{Appendix}
\renewcommand{\thesubsection}{\Alph{section}.\arabic{subsection}.}

\section{Derivation of (\ref{siki3-4})}
We present the derivation of (\ref{siki3-4}).
Assume that random variable pair $(x,y)$ follows a two-dimensional normal distribution and that each marginal distribution is a standard normal distribution.
We denote the correlation coefficient between $x$ and $y$ as $\rho$ and assume that $-1<\rho<1$.
Then, the distribution of $(x,y)$ is specified and the probability density function $p_{0,\rho}$ is described as
\begin{align}
p_{0,\rho}(x,y)=\frac{1}{2\pi\sqrt{1-\rho^2}}\exp\left(-\frac{x^2-2\rho xy+y^2}{2(1-\rho^2)}\right).
\end{align}
Assume that $(x_1,x_2),\cdots,(x_k,y_k)$ independently follow $p_{0,\rho}$.
For arbitrary positive integer $k$, we define $(X_k,Y_k)$ as
\begin{align}
X_k:=\sum_{j=1}^kx_j^2,\ Y_k:=\sum_{j=1}^ky_j^2.
\end{align}
We represent the probability density function of $(X_k,Y_k)$ as $p_{k,\rho}$, the characteristic function of $p_{k,\rho}$ as $\phi_{k,\rho}$ and the joint cumulative distribution function of $(X_k,Y_k)$ as $F_{k,\rho}$.
For arbitrary $k$, the probability that $X_k\leq0$ or $Y_k\leq0$ is zero, so in the following we consider the region where $X_k>0$ and $Y_k>0$.

By definition, $p_{1,\rho}$ satisfies
\begin{align}
\begin{split}
p_{1,\rho}(X_1,Y_1)
=&\frac{1}{4\sqrt{X_1Y_1}}\left\{p_{0,\rho}\left(\sqrt{X_1},\sqrt{Y_1}\right)+p_{0,\rho}\left(-\sqrt{X_1},\sqrt{Y_1}\right)\right.\\
&\quad\quad\quad\quad\quad\quad\quad\quad\left.+p_{0,\rho}\left(\sqrt{X_1},-\sqrt{Y_1}\right)+p_{0,\rho}\left(-\sqrt{X_1},-\sqrt{Y_1}\right)\right\}.
\end{split}
\end{align}
Then, the corresponding characteristic function $\phi_{1,\rho}$ is computed as follows:
\begin{align}
\phi_{1,\rho}(s,t)=&\int \int_{(0,\infty)\times(0,\infty)}p_{1,\rho}(X_1,Y_1)e^{isX_1+itY_1}dX_1dY_1\\
\begin{split}
=&\int \int_{(0,\infty)\times(0,\infty)}\frac{1}{4\sqrt{X_1Y_1}}\left\{p_{0,\rho}\left(\sqrt{X_1},\sqrt{Y_1}\right)+p_{0,\rho}\left(-\sqrt{X_1},\sqrt{Y_1}\right)\right.\\
&\quad\quad\quad\quad\left.+p_{0,\rho}\left(\sqrt{X_1},-\sqrt{Y_1}\right)+p_{0,\rho}\left(-\sqrt{X_1},-\sqrt{Y_1}\right)\right\}e^{isX_1+itY_1}dX_1dY_1.
\end{split}
\end{align}
Here, $i$ is the imaginary unit.
We introduce the following change of variables:
\begin{align}
\bar{X}=\sqrt{X_1},\ \bar{Y}=\sqrt{Y_1}.
\end{align}
Using $d\bar{X}=\frac{1}{2\sqrt{X_1}}dX_1$ and $d\bar{Y_1}=\frac{1}{2\sqrt{Y_1}}dY_1$, we get
\begin{align}
\begin{split}
\phi_{1,\rho}(s,t)=&\int \int_{(0,\infty)\times(0,\infty)}\Bigg\{p_{0,\rho}\left(\bar{X},\bar{Y}\right)+p_{0,\rho}\left(-\bar{X},\bar{Y}\right)\\
&\quad\quad\quad\quad+p_{0,\rho}\left(\bar{X},-\bar{Y}\right)+p_{0,\rho}\left(-\bar{X},-\bar{Y}\right)\Bigg\}e^{is\bar{X}^2+it\bar{Y}^2}d\bar{X}d\bar{Y}
\end{split}\\
=&\int \int_{(-\infty,\infty)\times(-\infty,\infty)}p_{0,\rho}\left(\bar{X},\bar{Y}\right)e^{is\bar{X}^2+it\bar{Y}^2}d\bar{X}d\bar{Y}\\
=&\int \int_{(-\infty,\infty)\times(-\infty,\infty)}\frac{1}{2\pi\sqrt{1-\rho^2}}\exp\left(-\frac{\bar{X}^2-2\rho \bar{X}\bar{Y}+\bar{Y}^2}{2(1-\rho^2)}\right)e^{is\bar{X}^2+it\bar{Y}^2}d\bar{X}d\bar{Y}\\
=&\frac{1}{2\pi\sqrt{1-\rho^2}}\int_{-\infty}^{\infty} \int_{-\infty}^{\infty}\exp\left(-\frac{\bar{X}^2-2\rho \bar{X}\bar{Y}+\bar{Y}^2}{2(1-\rho^2)}\right)e^{is\bar{X}^2+it\bar{Y}^2}d\bar{X}d\bar{Y}\\
\begin{split}
=&\frac{1}{2\pi\sqrt{1-\rho^2}}\int_{-\infty}^\infty \int_{-\infty}^\infty\exp\left(-\frac{1-2is(1-\rho^2)}{2(1-\rho^2)}\left\{\bar{X}-\frac{\rho\bar{Y}}{1-2is(1-\rho^2)}\right\}^2\right.\\
&\quad\quad \quad\quad \quad\quad\quad\quad  \left.-\frac{\left(1-2is(1-\rho^2)\right)\left(1-2it(1-\rho^2)\right)-\rho^2}{2(1-\rho^2)\left(1-2is(1-\rho^2)\right)}\bar{Y}^2\right)d\bar{X}d\bar{Y}
\end{split}\\
\begin{split}
\label{eq:1}
=&\frac{1}{2\pi\sqrt{1-\rho^2}}\int_{-\infty}^\infty \left[\exp\left(-\frac{\left(1-2is(1-\rho^2)\right)\left(1-2it(1-\rho^2)\right)-\rho^2}{2(1-\rho^2)\left(1-2is(1-\rho^2)\right)}\bar{Y}^2\right)\right.\\
&\quad\quad\quad\quad\quad\left.\times\int_{-\infty}^\infty\exp\left(-\frac{1-2is(1-\rho^2)}{2(1-\rho^2)}\left\{\bar{X}-\frac{\rho\bar{Y}}{1-2is(1-\rho^2)}\right\}^2\right) d\bar{X}\right]d\bar{Y}.
\end{split}
\end{align}
Since $-1<\rho<1$, we have
\begin{align}
&\text{Re}\left(\frac{1-2is(1-\rho^2)}{2(1-\rho^2)}\right)=\frac{1}{2(1-\rho^2)}>0,\label{appendix-4}\\
&\text{Re}\left(\frac{\left(1-2is(1-\rho^2)\right)\left(1-2it(1-\rho^2)\right)-\rho^2}{2(1-\rho^2)\left(1-2is(1-\rho^2)\right)}\right)
=\frac{1+4s^2(1-\rho^2)}{2\left(1+4s^2(1-\rho^2)^2\right)}>0.\label{appendix-5}
\end{align}
By (\ref{appendix-4}) and (\ref{appendix-5}), 
\begin{align}
&\int_{-\infty}^\infty\exp\left(-\frac{1-2is(1-\rho^2)}{2(1-\rho^2)}\left\{\bar{X}-\frac{\rho\bar{Y}}{1-2is(1-\rho^2)}\right\}^2\right) d\bar{X}=\sqrt{\frac{2\pi(1-\rho^2)}{1-2is(1-\rho^2)}},\label{eq:2}\\
\begin{split}
\label{eq:3}
&\int_{-\infty}^\infty \exp\left(-\frac{\left(1-2is(1-\rho^2)\right)\left(1-2it(1-\rho^2)\right)-\rho^2}{2(1-\rho^2)\left(1-2is(1-\rho^2)\right)}\bar{Y}^2\right)d\bar{Y}\\
&\quad\quad\quad\quad\quad\quad\quad\quad\quad\quad\quad\quad\quad\quad\quad=\sqrt{\frac{2\pi(1-\rho^2)\left(1-2is(1-\rho^2)\right)}{\left(1-2is(1-\rho^2)\right)\left(1-2it(1-\rho^2)\right)-\rho^2}}.
\end{split}
\end{align}
Here, the root of a complex number is defined as
\begin{align}
\sqrt{\alpha}:=\sqrt{|\alpha|}e^{\frac{i}{2}\text{Arg}\alpha}.
\end{align}
Substituting (\ref{eq:2}) and (\ref{eq:3}) into (\ref{eq:1}), we arrive at
\begin{align}
\begin{split}
\phi_{1,\rho}(s,t)=&\frac{1}{2\pi\sqrt{1-\rho^2}}\sqrt{\frac{2\pi(1-\rho^2)}{1-2is(1-\rho^2)}}\\
&\quad\quad\times\int_{-\infty}^\infty\exp\left(-\frac{\left(1-2is(1-\rho^2)\right)\left(1-2it(1-\rho^2)\right)-\rho^2}{2(1-\rho^2)\left(1-2is(1-\rho^2)\right)}\bar{Y}^2\right) d\bar{Y}
\end{split}\\
=&\sqrt{\frac{1-\rho^2}{\left(1-2is(1-\rho^2)\right)\left(1-2it(1-\rho^2)\right)-\rho^2}}.
\end{align}

\begin{figure}[h]
\begin{center}
\includegraphics[width=5cm]{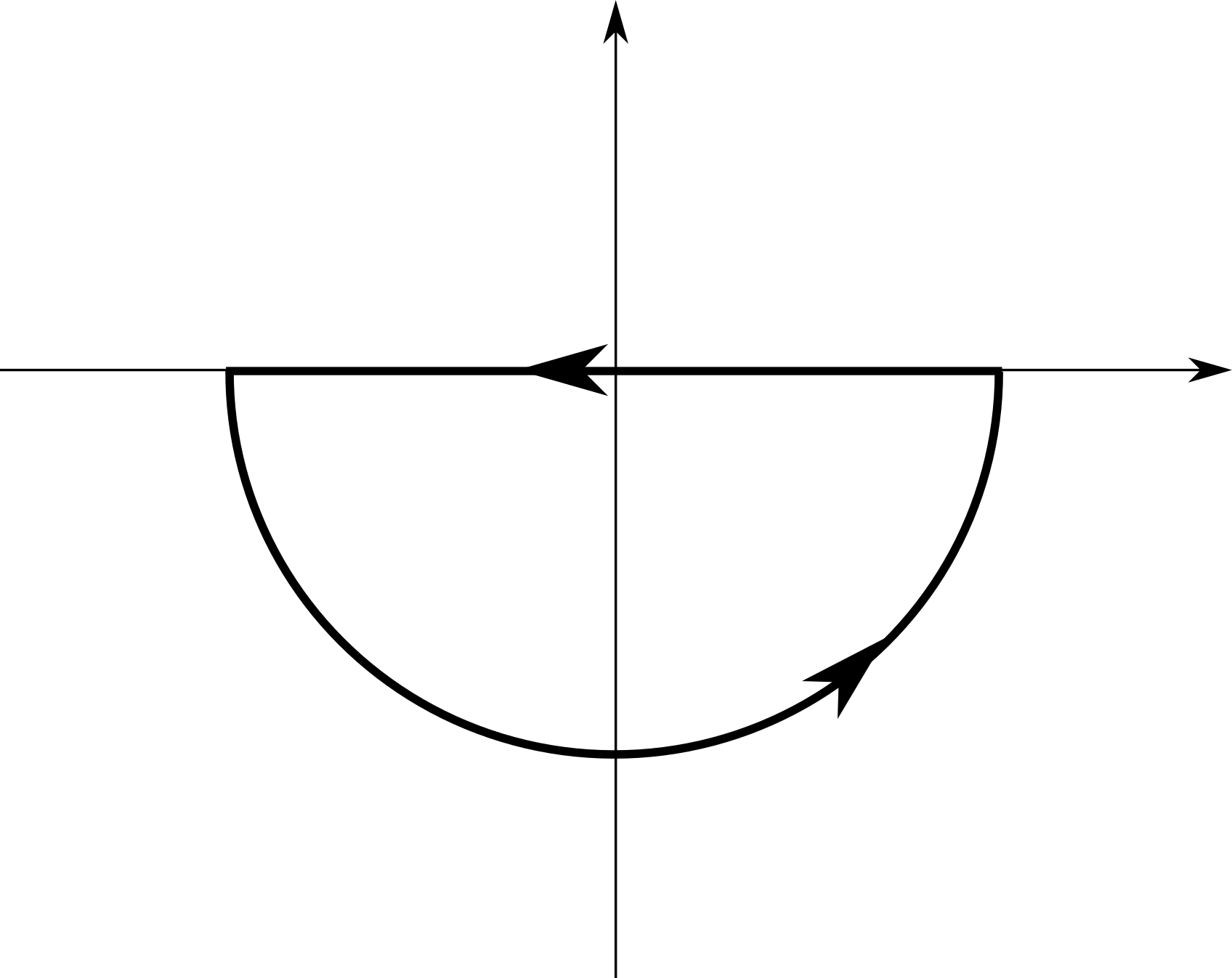}
\begin{picture}(0,0)
\put(-10,76){Re}
\put(-93,110){Im}
\put(-35,73){$R$}
\put(-130,73){$-R$}
\put(-135,30){$Re^{i\theta}$}
\end{picture}
\end{center}
\caption{Integral path $C_R$}
\label{path}
\end{figure}
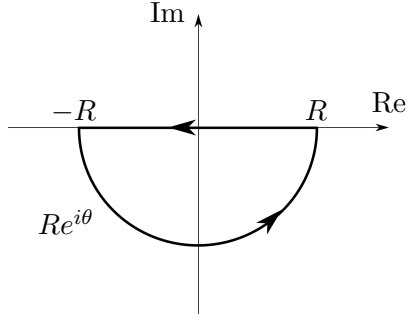

Since $(x_1,x_2),\cdots,(x_N,y_N)$ are mutually independent, we obtain
\begin{align}
\phi_{N,\rho}(s,t)=&\left[\phi_{1,\rho}(s,t)\right]^N\\
=&\frac{(1-\rho^2)^\frac{N}{2}}{\left\{\left(1-2is(1-\rho^2)\right)\left(1-2it(1-\rho^2)\right)-\rho^2\right\}^\frac{N}{2}}.
\end{align}
Then, we have
\begin{align}
p_{N,\rho}(X_N,Y_N)=&\frac{1}{4\pi^2}\int \int_{(-\infty,\infty)\times(-\infty,\infty)}\phi_{N,\rho}(s,t)e^{-isX_N-itY_N}dsdt\\
=&\frac{(1-\rho^2)^\frac{N}{2}}{4\pi^2}\int_{-\infty}^{\infty}\int_{-\infty}^{\infty}\frac{e^{-isX_N-itY_N}}{\left\{\left(1-2is(1-\rho^2)\right)\left(1-2it(1-\rho^2)\right)-\rho^2\right\}^\frac{N}{2}}dsdt\\
=&\frac{1}{4(-2i)^\frac{N}{2}\pi^2}\int_{-\infty}^{\infty}\frac{e^{-itY_N}}{\left(1-2it(1-\rho^2)\right)^{\frac{N}{2}}}\int_{-\infty}^{\infty}\frac{e^{-isX_N}}{\left\{s+\frac{i(1-2it)}{2\left(1-2it(1-\rho^2)\right)}\right\}^{\frac{N}{2}}}dsdt.
\end{align}
We assume $N$ to be even as per section 3 of the main text.
Since $-1<\rho<1$, 
\begin{align}
\text{Im}\left(-\frac{i(1-2it)}{2\left(1-2it(1-\rho^2)\right)}\right)=-\frac{1+4t^2(1-\rho^2)}{2\left(1+4t^2(1-\rho^2)^2\right)}<0.
\end{align}
Then，since $X_N>0$, using the integration path shown in fig. \ref{path}, Jordan's lemma can be applied as follows:
\begin{align}
&\int_{-\infty}^{\infty}\frac{e^{-isX_N}}{\left\{s+\frac{i(1-2it)}{2\left(1-2it(1-\rho^2)\right)}\right\}^\frac{N}{2}}ds\\
=&-\lim_{R\to\infty}\int_{C_R}\frac{e^{-isX_N}}{\left\{s+\frac{i(1-2it)}{2\left(1-2it(1-\rho^2)\right)}\right\}^\frac{N}{2}}ds\\
=&-\exp\left(-\frac{1-2it}{2\left(1-2it(1-\rho^2)\right)}X_N\right)\lim_{R\to\infty}\int_{C_R}\sum_{j=0}^\infty\frac{(-iX_N)^j}{j!}\left\{s+\frac{i(1-2it)}{2\left(1-2it(1-\rho^2)\right)}\right\}^{j-\frac{N}{2}}ds\\
=&\frac{2\pi (-i)^\frac{N}{2}{X_N}^{\frac{N}{2}-1}}{\Gamma\left(\frac{N}{2}\right)}\exp\left(-\frac{1-2it}{2\left(1-2it(1-\rho^2)\right)}X_N\right).
\end{align}
Then, we obtain
\begin{align}
p_{N,\rho}(X_N,Y_N)=&\frac{{X_N}^{\frac{N}{2}-1}}{2^{\frac{N}{2}+1}\pi\Gamma\left(\frac{N}{2}\right)}\int_{-\infty}^{\infty}\frac{\exp\left(-\frac{1-2it}{2\left(1-2it(1-\rho^2)\right)}X_N\right)e^{-itY_N}}{\left(1-2it(1-\rho^2)\right)^\frac{N}{2}}dt\\
=&\frac{{X_N}^{\frac{N}{2}-1}e^{-\frac{X_N}{2(1-\rho^2)}}}{(-i)^\frac{N}{2}2^{N+1}\pi\Gamma\left(\frac{N}{2}\right)(1-\rho^2)^\frac{N}{2}}\int_{-\infty}^{\infty}\frac{\exp\left(\frac{i\rho^2}{4(1-\rho^2)^2\left\{t+\frac{i}{2(1-\rho^2)}\right\}}X_N\right)e^{-itY_N}}{\left\{t+\frac{i}{2(1-\rho^2)}\right\}^\frac{N}{2}}dt.
\end{align}
Since $Y_N>0$, Jordan's lemma can be applied again as follows:
\begin{align}
&\int_{-\infty}^{\infty}\frac{\exp\left(\frac{i\rho^2}{4(1-\rho^2)^2\left\{t+\frac{i}{2(1-\rho^2)}\right\}}X_N\right)e^{-itY_N}}{\left\{t+\frac{i}{2(1-\rho^2)}\right\}^\frac{N}{2}}dt\\
=&-\lim_{R\to\infty}\int_{C_R}\frac{\exp\left(\frac{i\rho^2}{4(1-\rho^2)^2\left\{t+\frac{i}{2(1-\rho^2)}\right\}}X_N\right)e^{-itY_N}}{\left\{t+\frac{i}{2(1-\rho^2)}\right\}^\frac{N}{2}}dt\\
=&-e^{-\frac{Y_N}{2(1-\rho^2)}}\lim_{R\to\infty}\int_{C_R}\sum_{u=0}^\infty\sum_{v=0}^\infty\frac{1}{u!v!}\left(\frac{i\rho^2}{4(1-\rho^2)^2}X_N\right)^u(-iY_N)^v\left\{t+\frac{i}{2(1-\rho^2)}\right\}^{-u+v-\frac{N}{2}}dt\\
=&-2\pi ie^{-\frac{Y_N}{2(1-\rho^2)}}\sum_{r=0}^\infty\frac{1}{r!\Gamma\left(r+\frac{N}{2}\right)}\left(\frac{i\rho^2}{4(1-\rho^2)^2}X_N\right)^r(-iY_N)^{r+\frac{N}{2}-1}\\
=&2(-i)^\frac{N}{2}\pi {Y_N}^{\frac{N}{2}-1}e^{-\frac{Y_N}{2(1-\rho^2)}}\sum_{r=0}^\infty\frac{1}{r!\Gamma\left(r+\frac{N}{2}\right)}\left(\frac{\rho^2}{4(1-\rho^2)^2}X_NY_N\right)^r.
\end{align}
From the above, we arrive at
\begin{align}
p_{N,\rho}(X_N,Y_N)=&\frac{{X_N}^{\frac{N}{2}-1}{Y_N}^{\frac{N}{2}-1}}{2^k\Gamma\left(\frac{N}{2}\right)(1-\rho^2)^\frac{N}{2}}e^{-\frac{X_N+Y_N}{2(1-\rho^2)}}\sum_{r=0}^\infty\frac{1}{r!\Gamma\left(r+\frac{N}{2}\right)}\left(\frac{\rho^2}{4(1-\rho^2)^2}X_NY_N\right)^r.\label{appendix-6}
\end{align}

Finally, we derive $F_{N,\rho}$ from (\ref{appendix-6}).
For all non-negative integers $r$, we define the function $f_{r,\rho}$ as
\begin{align}
f_{r,\rho}(X_N,Y_N):=\frac{{X_N}^{\frac{N}{2}-1}{Y_N}^{\frac{N}{2}-1}}{2^N\Gamma\left(\frac{N}{2}\right)(1-\rho^2)^\frac{N}{2}}e^{-\frac{X_N+Y_N}{2(1-\rho^2)}}\frac{1}{r!\Gamma\left(r+\frac{N}{2}\right)}\left(\frac{\rho^2}{4(1-\rho^2)^2}X_NY_N\right)^r.
\end{align}
It is obvious that $f_{r,\rho}$ is continuous for all $r$.
In addition, for all $X_N>0$ and $Y_N>0$, $\sum_{r=0}^\infty f_{r,\rho}$ uniformly converges to $p_{N,\rho}$ on $(0,X_N]\times(0,Y_N]$.
Therefore, we can use termwise integration and calculate $F_{N,\rho}$ as follows:
\begin{align}
F_{N,\rho}(X_N,Y_N)=&\int \int_{(0,X_N]\times(0,Y_N]}p_{N,\rho}(X,Y)dXdY\\
=&\sum_{r=0}^\infty\int_0^{Y_N} \int_0^{X_N}f_{r,\rho}(X,Y)dXdY\\
\begin{split}
=&\sum_{r=0}^\infty\left[\frac{1}{2^N\Gamma\left(\frac{N}{2}\right)(1-\rho^2)^\frac{N}{2}r!\Gamma\left(r+\frac{N}{2}\right)}\left(\frac{\rho^2}{4(1-\rho^2)^2}\right)^r\right.\\
&\left.\quad\quad\quad\quad\times\int_0^{Y_N} Y^{r+\frac{N}{2}-1}e^{-\frac{Y}{2(1-\rho^2)}}dY\times\int_0^{X_N}X^{r+\frac{N}{2}-1}e^{-\frac{X}{2(1-\rho^2)}}dX\right].
\end{split}
\end{align}
We introduce the following notation: 
\begin{align}
X^\prime=\frac{X}{2(1-\rho^2)}.
\end{align}
Then, we obtain
\begin{align}
\int_0^{X_N}X^{r+\frac{N}{2}-1}e^{-\frac{X}{2(1-\rho^2)}}dX=&2^{r+\frac{N}{2}}(1-\rho^2)^{r+\frac{N}{2}}\int_0^{\frac{X_N}{2(1-\rho^2)}}{X^\prime}^{r+\frac{N}{2}-1}e^{-X^\prime}dX^\prime\\
=&2^{r+\frac{N}{2}}(1-\rho^2)^{r+\frac{N}{2}}\gamma\left(r+\frac{N}{2},\frac{X_N}{2(1-\rho^2)}\right),
\end{align}
where $\gamma$ is the incomplete gamma function. 
The integral on $Y$ can be computed similarly.
From the above, we get (\ref{siki3-4}).
\end{document}